\documentclass[a4,12pt]{amsart}
\usepackage{amsxtra,amsmath,amsthm,amssymb,amscd,amsfonts,
multicol,tabularx,latexsym}
\theoremstyle{definition}

\newtheorem{thm}{Theorem}
\newtheorem{lem}{Lemma}
\newtheorem{prop}{Proposition}
\newtheorem{rem}{Remark}

\newcommand{\Gal}{\mathrm{Gal}}

\newcommand{\Q}{\Bbb Q}
\newcommand{\Z}{\Bbb Z}

\newdimen\minCDarrowwidth
\date{}
\begin{document}
\title[]{Remark on the Alexander polynomials of
periodic knots}
\author[Manabu Ozaki]
{
Manabu\ OZAKI
}
\dedicatory{To Yuko's 29th birthday}
\maketitle
\begin{center}
Department of Mathematics,\\
School of Science and Engineering,\\
Kinki University,\\
Kowakae 3-4-1, Higashi-Osaka, 577-8502, JAPAN\\
e-mail:\ \verb+ozaki@math.kindai.ac.jp+
\end{center}
\section{Introduction}
Let $p$ be a prime number and $K$ a non-trivial knot
in $S^3$ which has period $p$.
Then the Alexander polynomial $\Delta_K(t)$ of such a knot
must have some distinguished properties as Murasugi \cite{Mur}
have revealed.

In this note, we shall show that under a certain assumption on
the Alexander polynomial $\Delta_K(t)$,
it is uniquely determined only by $p$.

The proof will be done by applying number theory
to Murasugi's work \cite{Mur} on a necessary condition on $\Delta_K(t)$
for periodic knots $K$.

\section{Result}
For a knot $K$
in $S^3$,
we denote by $\Delta_K(t)\in\Z[t]$ the Alexander polynomial
of $K$ normalized such that $\Delta_K(0)\ne 0$ and the leading
coefficient of it is positive.

Our main result is the following:
\begin{thm}\label{thm1}
Let $p$ be an odd prime number and $K\subseteq S^3$
a non-trivial periodic knot of period $p$.
If $\Delta_K(t)$ is monic and has degree $p-1$,
then we have
\[
\Delta_K(t)=\sum_{n=0}^{p-1}(-1)^nt^n=t^{p-1}-t^{p-2}+\cdots-t+1.
\] 
\end{thm}
For the case where degree $p-1$ Alexander polynomials of
$p$-periodic knots have general leading coefficients,
we will show the following:
\begin{thm}\label{thm2}
Let $p$ be an odd prime number and define
$\Pi(p)$ to be the set of all the $p$-periodic knots in $S^3$
whose Alexander polynomial has degree $p-1$.
Also, for a finite set $S$ of prime numbers which does not contain $p$,
we define the set $\mathcal{D}(p,S)\subseteq\Z[t]$ to be the
collection of all the $\Delta_K(t)$'s
such that $K\in\Pi(p)$ and the leading coefficient of
$\Delta_K(t)$ is prime to the prime numbers outside $S$.
Then $\mathcal{D}(p,S)$ is finite and
\[
\#\mathcal{D}(p,S)\le\left(\frac{p+1}{2}\right)
^{3\cdot 7^{\frac{3(p-1)}{2}+\#S(\Q(\zeta_p))}},
\]
where $S(\Q(\zeta_p))$ denotes the set of all the primes of
the $p$-th cyclotomic field $\Q(\zeta_p)$
lying over the prime numbers in $S$.

\end{thm}
\begin{rem}
(1)\ \ If $K$ is fibred, then $\Delta_K(t)$ is monic.

\noindent
(2)\ \ 
For a non-trivial knot $K$ of prime period $p$,
if the leading coefficient of $\Delta_K(t)$ is
prime to $p$, then $\deg \Delta_K(t)\ge p-1$
(See Davis and Livingston\cite[Cor. 4.2]{DL}).
\end{rem}

\section{Proof of Theorem 1.}
Let $T$ be a transformation of $S^3$ of order $p$ such that
$T(K)=K$ and it acts on $K$ fixed point freely, and we denote by
$B\subseteq S^3$ the set of the fix points of $T$.
Then $B$ is the unknot and the quotient space $S^3/T$ is homeomorphic to $S^3$,
and we let $\overline{K}$ and $\overline{B}$
be the quotient knots of $K$ and $B$ in $S^3/T$,
respectively.
We write for  $\Delta_{\overline{K}\cup\overline{B}}(t,u)\in
\Z[t,u]$
the two-variable Alexander polynomial of
the link $\overline{K}\cup\overline{B}$
with $\Delta_{\overline{K}\cup\overline{B}}(t,u)\not\in
t\Z[t,u]\cup u\Z[t,u]$
(Note that $\Delta_{\overline{K}\cup\overline{B}}(t,u)$
is defined up to $\pm 1$).
We put $\lambda$ the linking number of $K$ and $B$.

It follows from Murasugi \cite{Mur} 
that
\begin{equation}\label{dec}
\Delta_K(t)
=\Delta_{\overline{K}}(t)
\prod_{i=1}^{p-1}\Delta_{\overline{K}\cup\overline{B}}(t,\zeta_p^i)
\end{equation}
with a primitive $p$-th root of unity
$\zeta_p$.
Also, by using Murasugi's congruence \cite{Mur},
\[
\Delta_K(t)\equiv
\pm t^j\left(\frac{t^\lambda-1}{t-1}\right)^{p-1}\Delta_{\overline{K}}(t)^p
\pmod{p}
\]
for some $j\in\Z$,
we derive
$\lambda=2$ and $\Delta_{\overline{K}}(t)=1$,
because $\Delta_{\overline{K}}(t)\mid\Delta_{K}(t)$
in $\Z[t]$ by \eqref{dec}, 
$\deg\Delta_{K}(t)=p-1$, and the leading coefficient
of $\Delta_K(t)$ is prime to $p$.
Hence we may assume that
\begin{equation}\label{al}
\Delta_{\overline{K}\cup\overline{B}}(t,\zeta_p^i)=
g(\zeta_p^i)t-h(\zeta_p^i)
\end{equation}
with some $g(u), h(u)\in\Z[u]$ for $1\le i\le p-1$
by \eqref{dec}.

Because $\Delta_K(t)$ is monic, we see that $\eta_1:=g(\zeta_p)$
is a unit of the ring $\Z[\zeta_p]$ by 
the relation $\prod_{i=1}^{p-1}g(\zeta_p^i)=1$,
which comes from \eqref{dec}.
Also, since the constant term of $\Delta_K(t)$ is equal to 1,
we find that $\eta_2:=h(\zeta_p)$ is a unit of $\Z[\zeta_p]$. 
Hence if we put $\varepsilon:=\eta_1^{-1}\eta_2\in\Z[\zeta_p]^{\times}$,
then we have
\begin{equation}\label{eps}
\Delta_K(t)
=\prod_{i=1}^{p-1}g(\zeta_p^i)
\prod_{i=1}^{p-1}(t-g(\zeta_p^i)^{-1}h(\zeta_p^i))
=\prod_{\sigma\in\Gal(\Q(\zeta_p)/\Q)}(t-\sigma(\varepsilon)).
\end{equation}
On the other hand, it follows from the second Torres condition
(See \cite{Tor}) that
\begin{equation}\label{tor2}
\Delta_{\overline{K}\cup\overline{B}}(t^{-1},\zeta_p^{-1})
=t^a\zeta_p^b\Delta_{\overline{K}\cup\overline{B}}(t,\zeta_p)
\end{equation}
for some $a,b\in\Z$.
Then we find from \eqref{al} that
$a=-1$ and
\[
h(\zeta_p^{-1})=-\zeta_p^{b}g(\zeta_p),\ 
g(\zeta_p^{-1})=-\zeta_p^{b}h(\zeta_p). 
\]
Therefore, if we denote by $J\in\Gal(\Q(\zeta_p)/\Q)$
the complex conjugation, we have
\begin{equation}\label{J}
J(\varepsilon)=J(h(\zeta_p))/J(g(\zeta_p))
=h(\zeta_p^{-1})/g(\zeta_p^{-1})
=g(\zeta_p)/h(\zeta_p)=\varepsilon^{-1}.
\end{equation}
We need the following fact from the theory of cyclotomic fields
(see \cite[Prop.1.5]{Was}):
\begin{lem}\label{lem1}
For any $\varepsilon\in\Z[\zeta_p]^\times$, there exist
$r\in\Z$ and $\varepsilon_0\in\Z[\zeta_p+\zeta_p^{-1}]^\times$
such that
\[
\varepsilon=\zeta_p^r\varepsilon_0.
\]
\end{lem}
By Lemma \ref{lem1} and \eqref{J}, 
we obtain
\[
\varepsilon^{-1}=J(\varepsilon)
=
J(\zeta_p)^rJ(\varepsilon_0)
=
\zeta_p^{-r}\varepsilon_0=\zeta_p^{-2r}\varepsilon,
\]
from which we derive
\[
\varepsilon=\pm\zeta_p^r.
\]
Because $\Delta_K(1)=\pm 1$ and $\Delta_K(-1)\ne 0$,
we conclude that $\varepsilon=-\zeta_p^r$
with some $r\in\Z$ prime to $p$, and
\[
\Delta_K(t)=t^{p-1}-t^{p-2}+\cdots-t+1
\]
by \eqref{eps}. Thus we have proved Theorem 1.
\hfill$\Box$
\section{Proof of Theorem 2}

We will give the following proposition,
from which we can easily derive Theorem 2:
\begin{prop}
For a finite extension field $F/\Q$,
a finite set $S$ of prime numbers,
and positive integer $m$,
we define the set $\mathcal{P}(F,S,m)\subseteq\Z[t]$ to be
the collection of $\Delta_K(t)$'s for the knots $K$ in $S^3$
such that the leading coefficient of $\Delta_K(t)$
is prime to the prime numbers outside $S$,
the splitting field $\mathrm{Spl}(\Delta_K(t))$
of $\Delta_K(t)$ over $\Q$ is contained in $F$,
and the multiplicity of each zero of $\Delta_K(t)$
is at most $m$.
Then $\mathcal{P}(F,S,m)$ is finite
and we have 
\[
\#\mathcal{P}(F,S,m)\le (m+1)^{3\cdot 7^{[F:K]+\#(S(F)\cup\infty(F))}},
\]
where $S(F)$ and $\infty(F)$ denote
the set of the primes of $F$ lying over the prime numbers
in $S$ and that of the archimedean primes of $F$, respectively.
\end{prop}
{\bf Proof.}\ \ \ 
Let $\Delta_K(t)\in\mathcal{P}(F,S,m)$.
Then we have
\[
\Delta_K(t)=a\prod_{i=1}^d(X-\alpha_i)^{m_i}\in\Z[t]
\]
for some $a\in\Z$ which is prime to the prime numbers
outside $S$, $0\le d\in \Z$, distinct $\alpha_i$'s in $F$, and $1\le m_i\le m$.

Let $\frak{p}\not\in S(F)$ be any non-archimedean prime of $F$
and $v_\frak{p}$
a $\frak{p}$-adic valuation of $F$.
Since $a\alpha_i$ is integral over $\Z$ and $v_{\frak{p}}(a)=0$,
we see 
\[
v_\frak{p}(\alpha_i)=v_\frak{p}(a\alpha_i)\ge 0.
\]
Because the constant term and the leading coefficient
of $\Delta_K(t)$ are coincide,
we have
\[
a=\pm\prod_{i=1}^d\alpha_i^{m_i}, 
\]
from which we derive
\[
0=v_{\frak{p}}(a)=\sum_{i=1}^dm_iv_{\frak{p}}(\alpha_i). 
\]
Therefore we find that $v_{\frak{p}}(\alpha_i)=0$
for $1\le i\le d$, which means that $\alpha_i$'s are $S(F)$-units
of $F$.

On the other hand, 
$a(1-\alpha_i)=a-a\alpha_i$ is also integral over $\Z$,
we obtain
\[
v_{\frak{p}}(1-\alpha_i)=v_{\frak{p}}(a(1-\alpha_i))\ge 0.
\]
Also, since $\Delta_K(1)=\pm 1$,
we have
\[
a\prod_{i=1}^d(1-\alpha_i)^{m_i}=\pm 1,
\]
from which we derive 
\[
0=v_{\frak{p}}(\pm 1)=v_{\frak{p}}(a)+\sum_{i=1}^dm_iv_{\frak{p}}(1-\alpha_i)
=\sum_{i=1}^dm_iv_{\frak{p}}(1-\alpha_i).
\]
Hence $v_{\frak{p}}(1-\alpha_i)=0$ for $1\le i\le d$,
 which means that $1-\alpha_i$'s are also $S(F)$-units
of $F$.

Now we apply the following result from analytic number theory
given by Evertse \cite{Eve}:
\begin{lem}\label{lem2}
Let $F$ be a finite extension of $\Q$ and $T$ a finite
set of non-archimedean primes of $F$.
Then the number of solutions $(X,Y)$ of the equation
\[
X+Y=1
\]
in the $T$-unit group of $F$ is at most
$3\cdot 7^{[F:\Q]+\#(T\cup\infty_F)}$.
\end{lem}
As we have seen in the above, $(\alpha_i, 1-\alpha_i)$
is a solution in the $S(F)$-unit group of $F$
of the equation $X+Y=1$.
Hence, it follows from Lemma \ref{lem2}
that the number of such $\alpha_i$'s
is at most $3\cdot 7^{[F:\Q]+\#(T\cup\infty_F)}$.
Therefore we obtain
\[
\#\mathcal{P}(F,S,m)\le (m+1)^{3\cdot 7^{[F:\Q]+\#(T\cup\infty_F)}}.
\]
\hfill$\Box$

Now we will derive Theorem 2 from Proposition 1.
Assume $K\in\mathcal{D}(p,S)$.
Then, as the proof of Theorem 1, we find that
\begin{equation*}
\Delta_K(t)=\prod_{i=1}^{p-1}(g(\zeta_p)t-h(\zeta_p))
\end{equation*}
for some $g(u), h(u)\in \Z[u]$,
since the leading coefficient of $\Delta_K(t)$ is
prime to $p$ by $p\not\in S$ and $\deg\Delta_K(t)=p-1$.
Hence $\mathrm{Spl}(\Delta_K(t))\subseteq\Q(\zeta_p)$
and $\Delta_K(t)\in\mathcal{P}(\Q(\zeta_p),S,\frac{p-1}{2})$
because $\Delta_K(t)$ has at least two distinct zeros.
Therefore, applying Proposition 1,
we complete the proof of Theorem 2 by using
the facts $[\Q(\zeta_p):\Q]=p-1$ and $\#\infty_{\Q(\zeta_p)}
=\frac{p-1}{2}$.
\hfill$\Box$

\end{document}